# A K-Theory approach to characterize admissible physical manifolds

Patrick Linker (1), Cenap Ozel (2), Alexander Pigazzini (3), Monika Sati (4), Richard Pincak (5) and Eric Choi (6)

(1) Fachbereich Kontinuumsmechanik, Technische Universität Darmstadt, 64289 Darmstadt, Germany;
    E-Mail: mrpatricklinker@gmail.com
(2) Department of Mathematics, King Abdulaziz University, Jeddah 21589, Saudi Arabia;
    E-Mail: cozel@kau.edu.sa
(3) Mathematical and Physical Science Foundation, 4200 Slagelse, Denmark;
    E-Mail: pigazzini@topositus.com
(4) Department of Mathematics, Govt PG College, Joshimath, Chamoli, Uttarakhand;
    E-Mail: monikasati123@gmail.com
(5) Institute of Experimental Physics, Slovak Academy of Sciences, Kosie, Slovak Republic
    E-Mail: pincak@saske.sk
(6) School of Science and Technology, Georgia Georgia Gwinnett College, 1000 University Centre LN, Lawrenceville, GA 30043, USA
    E-Mail: ericchoi314@gmail.com

## Abstract

We will classify physically admissible manifold structures by the use of Waldhausen categories. These categories give rise to algebraic K-Theory. Moreover, we will show that a universal K-spectrum is necessary for a physical manifold being admissible. Application to the generalized structure of D-branes are also provided. This might give novel insights in how the manifold structure in String and M-Theory looks like.

**Math Subject Classification:** 18F25, 19M05, 55U15, 55R35, 81T30.

## Introduction

The Standard model has proven to successfully describe fundamental interactions of nature. However, several open questions in particle physics arise, concerning phenomena that cannot be described by the Standard Model [1]. One of these open questions are how to unify General Relativity, the classical theory of gravitation with Quantum Field theory. Such a unification is called Quantum Gravity.

Several models for the quantization of gravity were proposed. One of these approaches is String theory, that assumes particles not to be pointlike, but one-dimensional extended objects called strings. These strings can vibrate and depending on their vibration, these represent some sort of elementary particle. It assumes that the spacetime has more dimensions than the usual 3+1 dimensions of Lorentzian spacetime; however, these extra dimensions appear to be compactified in the limiting case of the Standard model. Compactification of higher-dimensional spacetime can, under certain circumstances [2], resemble fundamental forces of nature, for example in Kaluza-Klein theory, electromagnetism appears to be gravity in a 5-dimensional spacetime in which the fifth dimension is compactified on a circle. Another approach to Quantum Gravity is Loop quantum gravity [3]. This approach uses spin networks and its evolution to quantize space.

We will reformulate ordinary Quantum Field Theory on a continuous base manifold as a discrete manifold by introducing discrete chunks of the spacetimes, which can be represented by some more compact topology. An example of such a compact topology would be the topology of multiple compact groups.

There are Quantum Field Theories called Group Field Theory (GFT), meaning that its base manifold will be a compact group instead of a continuous manifold. There are approaches to Quantum Gravity in terms of Group Field Theory [4]. By choosing a representation basis for the groups, one can get a generalized matrix model. After treating the mathematics of the transition to a Quantum Field Theory with compact topology we will provide examples like how one can reformulate an ordinary Quantum Field Theory in terms of a Group Field Theory. Our mathematical construction will be a universal K-theory. This K-theory will be applied to D-branes.

## 1: Construction of the new field manifold

For a physical theory to be meaningful, we have to define the Hilbert space corresponding to all corresponding configurations $\mathcal{H}$. Moreover, we will have a couple of symmetries that we describe with a set of symmetry group actions $\mathcal{g}$ that acts on this Hilbert space. Let the manifold on which the fields are defined be $M$. In this paper we will focus on the ordinary 3+1 (curved) spacetime which is locally isomorphic to the Lorentzian spacetime $M \cong \mathbb{R}^4$. It is supposed to be a smooth $C^\infty$ manifold. Fields on a point in spacetime can be regarded as a section $\Gamma$ of a fiber bundle $E$ with projection map $\pi\colon E \to M$ and whose fibers are the vector space $W$, dependent on the type of the quantum field. Because there is some symmetry group action for example $U(1)$ for electromagnetic fields, we will focus on principal bundles with group actions $g\colon \mathcal{g} \times E \to E$.

We will foliate the principal bundle defined above into a new fiber bundle. The foliation will map into a fully discrete space which is locally isomorphic to $I \times R$, where $I$ is the discrete space on which a group action with elements in a finite group $\mathcal{g}^*$ exist and $R$ is the residual discrete space with $\#R = k$ and $k$ the number of indices the matrix-valued field has.

To achieve the transition from a smooth bundle to a discrete bundle, we will introduce the exact category, which is defined as follows:

**Definition:** An exact category $\mathcal{E}$, a subcategory of an abelian category $Ab$ is an additive category with objects $L', L, L''$ and morphisms $f_1, f_2$ such that there exists the short exact sequence

$$0 \to L' \xrightarrow{f_1} L \xrightarrow{f_2} L'' \to 0 \ (1a)$$

and it contains the split exact sequence

$$0 \to L' \xrightarrow{f_1} L' \oplus L'' \xrightarrow{f_2} L'' \to 0. \ (1b)$$

This category is closed under isomorphisms and closed under extensions. $f_1$ is called an admissible monomorphism and $f_2$ is called an admissible epimorphism. We will also denote exact category with a pair $(\mathcal{E}, \Sigma)$, where $\Sigma$ is the class of short exact sequences.

We will associate the category of fiber bundles with another kind of category called Waldhausen category that is defined as follows:

**Definition:** The Waldhausen category is the triple $(C, co(C), we(C))$ with a category $C$, the morphism $co(C)$ which is called cofibration and and a morphism $we(C)$ which is called weak equivalence such that $C$ contains the zero element, both $co(C), we(C)$ contain isomorphisms, are closed under compositions and are compatible with pushouts. Moreover, the unique cofibration $c: 0 \mapsto A \in co(C)$ for each object $A \in C$ exists.

Let $X$ be a contractible topological space (including fiber bundles) such that we can homeomorphically contract it to a point $\{*\}$. Then we can promote this topological space as an object of a Waldhausen category $(C_t, co(C_t), we(C_t))$ associated with category of topological spaces $C_t$, by defining $q: X \to \{*\} \in we(C_t)$, the contraction homeomorphism as a weak equivalence. A cofibration $\eta \in co(C_t)$ is defined to add some group structure on it. In a topological space without group structure, there exist morphisms between topological spaces $f$ and also, in topological spaces that contain group structure, we have the morphisms $f_g$. We then can introduce a cofibration $i$ that satisfies $f \circ i = f_g$. Denote $X_g$ a topological space that is cofibrated with a group structure $g$. Then the following pushout holds:

$$\begin{array}{ccc} X & \xrightarrow{\eta} & X_g \\ \downarrow q & & \downarrow q' \\ \{*\} & \xrightarrow{\eta} & X_g \cup_X \{*\}. \end{array} \quad (1c)$$

Since $X_g \cup_X \{*\} \cong X_g$ we conclude that the weak equivalence $q' \in we(C)$ is an isomorphism, in terms of topological spaces, a homeomorphism. The inclusion of isomorphisms in $we(C)$ is provided; likewise if we choose $g = \{id\}$, the trivial group containing only the identity element, we would have an isomorphism in $co(C)$.

Exact categories, also called Quillen's exact categories, are related with Waldhausen categories in terms of Waldhausen categories generalize the exact category.

**Lemma:** All exact categories are Waldhausen categories.

The crucial property of Waldhausen categories is that they can be promoted by Waldhausen's S-construction to a K-Theory. This is achieved by introducing a functor $S_m$ for an integer $m \in \mathbb{N}$ that associates a $m$-simplex to the Waldhausen category $C$. Denote $\omega S_m C$ the $m$-simplicial subcategory of weak equivalences. Denote $N$ the nerve of a category. Then the $m$-th K spectrum is defined by

$$K(C)_m = |N\omega S_m C|, \quad (1d)$$

where || denotes a topological realization. We have now constructed an algebraic K-theory that is associated with both the continuous spacetime manifold bundle and the discretized version that is isomorphic to $I \times R$ of it.

Let $\Delta: E \to E'$ be some triangulation of the fiber bundle $E$ to a discrete fiber bundle $E'$ which is locally isomorphic to $I \times R$. According to (1c) both can be associated with a Waldhausen category. The S-construction will generate the same algebraic K-theory associated with these. As the K-theory associated with such kind of Waldhausen categories with topological spaces as objects can be quite general, we should do some restrictions on these topological spaces. Note that in Quantum Field Theory, the fields should be defined on a Hausdorff space equipped with a metric and a connection. We propose that a physical topological space contains the de-Rham complex $\Lambda_*$ or some generalization on it. Let $\Lambda_n \otimes \mathcal{R}$ be the space of $n$-forms expressed in some ring $\mathcal{R}$ (e.g. a Lie algebra), then we have a long exact sequence

$$0 \to \cdots \to^d \Lambda_n \otimes \mathcal{R} \to^d \Lambda_{n+1} \otimes \mathcal{R} \to^d \Lambda_{n+2} \otimes \mathcal{R} \to^d \ldots \to 0. \quad (1e)$$

Here, $d$ is the exterior covariant derivative satisfying $d^2$ (Bianchi identity). To any field variable $\phi \in \Lambda_n \otimes \mathcal{R}$ that appears in a quantum field theory, we can assign the action of a short exact sequence on it

$$0 \to \psi \to^d \phi \to^d d\phi \to 0, \quad (1f)$$

where $\psi$ is some potential field such that $d\psi = \phi$ that maps injectively. This potential field generates a generalized gauge invariance for the field $d\phi$. The identity (1f) gives rise to restrict the Waldhausen category with topological spaces as object to the exact category (1a) as the Waldhausen category where the K-theory can be constructed. We will project the algebraic K-theory defined by (1d) to these constructed from exact categories. Denote these K-theory by $K^*$ it holds analogously to (1d):

$$K^*(C)_m = |N\omega S_m \mathcal{E}|. \quad (1g)$$

In n-dimensional spacetime, we define actions on $\Lambda_n \otimes \mathcal{R}$. Moreover we have gauge symmetries like Lorentz invariance. These ingredients are encoded in (1e) and (1f), such that $\Lambda_n \otimes \mathcal{R} \in \mathcal{E}$. Thus, the physically relevant K-theory is $K^*$.

We discretize the spacetime and its physical action functional such that one gets chunks $\Delta M$ on it and replacing the coordinate charts of the atlas of base manifold $\Delta M$ with compact topological spaces. We will regard the resulting differential-geometric structure as a Waldhaus category by construction (1c) and construct the same K-theory (1g) with help of Waldhaus S-construction. Thus, the K-theory is preserved by this discretization.

**Theorem:** A discretization $\Delta: E \to E'$ of fiber bundles exist if and only if the reduced K-theory $K^*$ associated to $E$, defined by (1g), is preserved.

One performs a functor $\Pi: E \to K^*$ from fiber bundles to algebraic K-theories to get the K-spectrum of the ordinary smooth spacetime manifold and its fields on it. Due to (1f), gauge symmetries also affect the K-spectrum. The new discrete bundle $E'$ which can define an action that is a generalized matrix model (for example GFT) will also been sent to the K-theory spectrum by $\Pi: E' \to K^*$. Now $E'$ is a finite set and the K-theory spectrum is also finite. On the other hand, $\Pi: E \to K^*$ forgets about the continuity structure of $E$. Therefore, $\Pi$ is a forgetting functor with respect to a continuous structure. However, information of the metrical and continuous structure will not get lost since we mapped the fiber bundles to a restricted K-theory spectrum induced by exact categories.

The dimension of matrices $\#R$ appear to be the necessary number of degree of freedom such that its corresponding K-theory agrees with that of $E$, the continuous fiber bundle. The symmetry group dimension $I$ is chosen such that the underlying K-theory is based on exact category $\mathcal{E}$ instead of some more general category $C$.

An alternative way to prove the theorem of preserving the K-groups is the following: Suppose that $\mathcal{E}$ is the category of continuous fiber bundles where deRham cohomology can be defined, one can create a sequence of cofibrations depicted by arrows between objects $A_{ij} \in \mathcal{E}$ for the horizontal index $i \in \{1,2,\dots,n\}$ and a vertical index $j \in \{1,2,\dots,n\}$ if there is constructed the $n$-th S-construction $S_n\mathcal{E}$. The following diagram commutes:

$$* \to A_{11} \to A_{21} \to \cdots \to A_{n1}$$
$$\downarrow \quad \downarrow \quad \downarrow \quad \downarrow$$
$$* \to A_{22} \to \dots \to A_{n2}$$
$$\downarrow \quad \downarrow \quad \downarrow$$

$$\begin{array}{ccc} \ldots \rightarrow & \ldots \rightarrow & \cdots \\ \downarrow & \downarrow & \\ * & A_{nn} & \\ & \downarrow & \\ & * & \text{(1h)}. \end{array}$$

In a physical theory we have symmetry groups like Poincaré symmetry and several gauge symmetries that act on the spacetime manifold. Cofibrations can be performed by just adding some group actions which leave the physical theory invariant. In (1h), horizontal arrows represent group actions on a principal bundle, while vertical arrows give rise to weak equivalences. Since group actions pose invariances, we have all arrows as homotopies between spaces. The S-construction where every arrow is some kind of homotopy equivalence (this includes the weak equivalence) gives rise to the physical meaningful algebraic K-theory. Define a functor $\mathcal{M}: \mathcal{E} \rightarrow \mathcal{E}^*$ that preserves the homotopy and weak equivalences and the S-construction (1h). Also, the cohomology structure is preserved, such that same exact sequences in $\mathcal{E}^*$ can be defined. Then we preserve all physically relevant symmetries. Moreover, the algebraic K-theory that the S-construction induces, will be the same as in $\mathcal{E}$. This proves that two different forms of a quantum field theory must have the same algebraic K-theory.

As an example, we want to express a 3+1-dimensional Quantum Field Theory in terms of discrete chunks, where each chunk can be covered by several compact groups. By assigning the functor $\mathcal{M}$ as sending a deRham chain complex to a chain complex of Lie group cohomology chain complex, we conclude that the dimension of the coordinate basis in continuum field theory must be the same as the basis in a group field theory. Let $\phi_I$ be a quantum field in a bounded spacetime chunk with compact support labeled by index $I$. In continuum spacetime, the field $\phi$ would carry 4 arguments, 1 time and 3 space arguments. By functoriality, in the group field theory framework it must depend on 4 group arguments, meaning that $\phi_I = \phi_I(g_0, g_1, g_2, g_3)$ where $g_0, g_1, g_2, g_3 \in G$ for some compact Lie group $G$. We will express a field $\phi(x)$ on some spacetime point $x \in M$ in the form

$$\phi(x) = \sum_I \int dg_0 \int dg_1 \int dg_2 \int dg_3 \phi_I(g_0, g_1, g_2, g_3) \lambda_I(x; g_0, g_1, g_2, g_3). \quad (1i)$$

Here, $\int dg_i$ for index $i \in \{0,1,2,3\}$ is the Haar measure of the corresponding group and $\lambda_I(x; g_0, g_1, g_2, g_3)$ are the coefficients that express spacetime location $x$ in terms of the region indexed by $I$ and the four group elements $g_0, g_1, g_2, g_3$. These coefficients are clearly zero if $x$ is not contained in region $I$ and nonzero otherwise. Existence of the Haar invariant measure is guaranteed, since any group actions will leave the physical theory invariant; the image of the functor $\mathcal{M}$ will preserve invariances under group actions and weak equivalences. If (1i) is substituted into an action functional of a regular continuum Quantum Field Theory one will have continuous spacetime integrations only over polynomials in $\lambda_I$ and its derivatives which result in some coefficient independent on any physical field. The quantum field theory will then be a group field theory with $\#R = 4$, i.e. a 4-dimensional matrix model.

## 2: D-branes in the framework of an algebraic K-theory

D-branes are objects in String Theory on which Strings can be attached to restrict the degree of freedom that Strings can take. These are crucial in shaping the structure of the string mass spectrum and the general dynamics of it. Suppose that $\Omega_p$ is a $p$-brane, i.e. an extended $p$-dimensional object. We can define a cofibration $c: \Omega_p \to \Omega_p \times \Sigma$, where $\Sigma$ is the open string worldsheet. All diffeomorphisms in $\Omega_p \times \Sigma$ define further cofibrations. Denote $\Pi_{ij}$ a construction isomorphic to $\Omega_p \times \Sigma$, where $i$ is the $i$-th configuration in string worldsheet and $j$ is the $j$-th D-brane configuration. Analogous to the cofibration sequences (1h), we can formulate the commjuting triangle of cofibration sequence

$$* \to \Pi_{11} \to \Pi_{21} \to \cdots \to \Pi_{n1}$$
$$\downarrow \quad \downarrow \quad \downarrow \quad \downarrow$$
$$* \to \Pi_{22} \to \ldots \to \Pi_{n2}$$
$$\downarrow \quad \downarrow \quad \downarrow$$
$$\ldots \to \ldots \to \cdots$$
$$\downarrow \quad \downarrow$$
$$* \quad \Pi_{nn}$$
$$\downarrow$$
$$* \quad (2a).$$

Hence, we obtain the same algebraic K-theory from (2a) as that of ordinary Quantum Field theory. We have ruled out that string worldsheet deformations are represented by "horizontal arrow" cofibration morphisms, while dynamics of D-branes are represented by "vertical arrow" cofibration morphisms. D-branes just determine what kind of motions a string worldsheet can perform. They act as a categorical pushout generator as one can see from (2a).

The String and D-brane configuration that generates categorically a cofibration (2a) is a very general construction. We know that only a particular class of Quantum Field theories can make sense; these are those that possess some kind of symmetries like Poincaré invariance. Therefore, the construction (2a) might not be adequate to Quantum Field Theories that possess symmetries like Poincaré invariance and gauge symmetries. Likewise, the construction (2a) should be restricted to admissible String-Brane configurations.

From String Theory it is known that multiple branes, where each brane incorporates a $U(1)$ symmetry on which strings can attach in all possible combinations, resemble, depending on the configuration of Strings and Branes, gauge invariances with gauge groups that are subgroups of $U(N)$. Here, $N$ is the number of $U(1)$ branes, or more generally, a rank of the underlying Lie group. For example, if there are two D-branes with strings that can either attach one of these branes with both ends or strings that connect these two branes, one has a $U(2)$ gauge symmetry if both branes intersect and $U(1) \times U(1)$ symmetry otherwise. The mathematics behind this construction is the Chan-Paton bundle $\nabla_B$. It is defined as a map from spacetime $X$ to the circle 2-group

$$\nabla_B : X \to B^2 U(1), \qquad (2b)$$

Where $B$ denotes the delooping operator. Consider the Lie Group Extension

$$U(1) \to U(N) \to PU(N). \qquad (2c)$$

$PU(N)$ is the projective unitary group of dimension $N$. From (2c) we can construct the homotopy fiber sequence

$$U(1) \to U(N) \to PU(N) \to BU(1) \to BU(N)$$
$$\to BPU(N) \to^{d_N} B^2 U(1). \qquad (2d)$$

In (2d), the map $d_N$ is the smooth refinement of the Dixmier-Douady class ($d_N : BPU(N) \to K(\mathbb{Z}, 3)$ with Eilenberg-MacLane space $K$). Now, let $c : X \to B^2 U(1)$ be a principal 2-group-2-bundle and $\sigma : c \to d_N$ a map to the smooth refinement. Then, the following diagram commutes:

$$X \to BPU(N)$$
$$\searrow_c \qquad \downarrow_{d_N}$$
$$B^2 U(1)$$

(2e).

There is also an isomorphism $BPU(N) \cong (BU(N)//BU(1))$ where $(A//B)$ denotes the categorical quotient between $A$ and $B$. Now we define the brane immersed in spacetime with the brane region $Q$ as the inclusion map $\iota: Q \to X$ and hence we get the following commutative diagram:

$$\begin{array}{ccc} Q & \xrightarrow{\nabla_g} & (BU(N)//BU(1)) \\ \downarrow_\iota & & \downarrow_{d_N} \\ X & \xrightarrow{\nabla_B} & B^2 U(1) \end{array}$$

(2f).

The gauge bundle map is denoted by $\nabla_g$. This is now a twisted bundle with connection on a D-brane with cohomology class $H^3(X, \mathbb{Z})$.

Now we do the S-construction with restriction to Chan-Paton bundles. Since we have defined twisted bundles on D-brane, we can associate these with vertical arrows of (2a). In (2f), the vertical arrows depict the mathematical structure of the D-branes that are governed by smooth refinements $d_N$ and inclusions in spacetime. The commutative diagram (2f) can be generalized to D-branes that also have strings attached to it. Open strings are also 1-branes, but without explicitly associating a gauge group on it. The endpoints of these strings, however, must lie in a D-brane and have the twisted bundle structure of it. Strings that connect D-branes can be regarded as a morphism between these D-branes if the D-branes are the object of some category $\mathcal{B}$ of D-branes. From previous considerations, string configurations (i.e. where strings begin and end) are depicted by horizontal arrows in (2a). Beginning points of strings are sources and endpoints of strings are the targets. Let the commutative diagrams in (2f) be objects of $\mathcal{B}$. Then morphisms in $\mathcal{B}$ are defined by

$$m_Q: Q \to Q', \quad (2g)$$

where $Q'$ is the brane region that obeys the same twisted bundle structure (2f) as $Q$. Therefore, $m_Q$ is a group action on $(BU(N)//BU(1))$ and $B^2U(1)$ that leaves the twisted bundle invariant; it preserves the twisted bundle structure. Horizontal arrows in the S-construction (2a) are morphisms of the form (2g). From the exact sequence (2d) we observe that $d_N$ is a weak equivalence if and only if $(BPU(N)//BU(N)) \cong BPU(N)$, meaning that groups acting via morphism (2g) must lie in $BU(N)$. To summarize, we can do the S-construction that also leads to the same algebraic K-theory as ordinary Quantum Field Theory if we define cofibration morphisms between twisted principal bundles with bundle-preserving morphisms (2g) and group actions on this twisted bundle in $BU(N)$. The loop space of $BU(N)$ is $U(N)$. Since we know that open strings $\Sigma$, that are some 2-dimensional topological spaces must carry group actions in $BU(N)$. If $\Omega$ maps some topological space to the corresponding loop space, $\nabla_g$ is a gauge bundle map and $\nabla_{g'}$ is a gauge bundle map in loop spaces, then the following diagram commutes:

$$\Sigma \xrightarrow{\nabla_g} BU(N)$$
$$\downarrow \Omega \qquad \downarrow \Omega$$
$$\Omega\Sigma \xrightarrow{\nabla_{g'}} U(N) \qquad (2h).$$

The loop space of the open string region $\Omega\Sigma$ is nontrivial, if the open string begins and ends on the same point. This is the case if and only if it begins and ends on the same D-brane or if it ends on another D-brane, that intersects the D-brane on which the string begins. From (2h), this string must have gauge group $U(N)$. Gauge groups of the Standard model like $U(1) \times SU(2)$ for electroweak interaction and $SU(3)$ for the strong interaction are subgroups of $U(N)$. We declare $\Omega\Sigma$ to be the string in lowest energy configuration, since it carries minimum energy if it is unstretched, i.e. a point. Thus, we expect that in low-energy regime, physical particles that satisfy the K-theory condition, must be subgroups of $U(N)$.

## 3: Possible discrete D-brane construction with PNDP-manifolds

As said, we want reformulate ordinary Quantum Field Theory on a continuous base manifold as a discrete manifold by introducing discrete chunks of the spacetimes.

For example, in [5], the authors consider a special manifold, called PNDP-manifold (Partial Negative Dimensional Product manifold), to try to describe the discrete gravity in emerging space theory, and also in this context we want to consider such manifolds as an example of discrete brane construction.

Since a PNDP-manifold (which is a special case of Einstein sequential warped-product manifold), has the fiber-manifold (which from the point of view of derived geometry and Kuranishi neighborhood) with negative "virtual" dimension *m* (therefore *m* will be an integer value negative), from a speculative point of view, also the dimension of the whole PNDP-manifold *(M,gM)* will be virtual (i.e. *dim(M) = dim(B)+dim(F) = dim(B)+m*, where we have *dim(B) = dim(B1)+dim(B2)*, because *B = B1 × B2*).

So, depending on the value of *m*, *dim(M)* could be "virtual" positive, zero, or negative. This means, again from a speculative point of view, that the "virtual" negative dimensions of F (fiber-manifold) interact and "virtually" cancel each other out with the positive dimensions of B (base-manifold).

We believe that space is not a fundamental property of dimensions. Space, we believe, may be a secondary property created by more fundamental entities and in this sense dimensions too could vanish via interactions with each other.

Now, if we consider a PNDP-manifold *M:= (I1 X I2) X [(I3 X I4) + E]*, where the fiber-manifold (*F:= [(I3 XI4) + E]*), is considered as a Kuranishi neighborhood *(I_3 X I_4, E, S)*, with product manifold *(I_3 X*

I_4), obstruction bundle $E \to (I_3 \times I_4)$, and section $S : (I_3 \times I_4) \to E$, then the dimension of the derived smooth fiber-manifold F is equal to $dim(I_3 \times I_4) - rank(E)$ (see [6] for more details about Kuranishi neighborhood and obstruction bundle), and where, in this specific case, $I_i \in R$ (with i=1, 2, 3, 4) are closed intervals with Planck length and $rank(E) = 4$. Then $dim(F) = -2$, so $dim(M) = 2 + (-2) = 0$, and thus, M emerges as a "virtual" point-like manifold.

We report the following well-known proposition:

**Proposition:** A topological space is a *0*-manifold if and only if it is a countable discrete space.

**Proof.** *Suppose that N is 0-manifold. Since it is a 0-manifold, for each $p \in N$, there must exist a neighborhood U that is homeomorphic to an open set in $R^0$. The only open sets in $R^0$ are $\emptyset$ and $\{0\}$. Since a nonempty set cannot be mapped into the empty set, U must be homeomorphic to $\{0\}$. If U contains more than 1 element, then any map from U into $\{0\}$ is not injective. Thus U must only contain p.*

*This shows that $\{p\}$ is open for every $p \in N$, so N has a discrete topology.*

*Let Q be a countable basis of N. Let $p \in p$ be given. Then $\{p\}$ is the union of elements of Q by definition. This implies that $\{p\} \in Q$.*

*Since Q contains all singletons and Q is countable, there are at most countably many singletons. In other words, N must be countable.*

*Therefore, N is a countable discrete space.*

*On the other hand, suppose that N is a countable discrete space.*

*$\{\{p\}|p \in N\}$ is a countable basis.*

*N is Hausdorff since any $p_i \neq p_j$ (with $i \neq j$), can be separated by $\{p_i\}$ and $\{p_j\}$.*

*For every point p ∈ N, {p} is a neighborhood of p that is locally homeomorphic to R^0. Thus N is a 0-manifold.* □

Considering a discrete 0-dimensional D-brane (N), composed of "virtual" point-like PNDP-manifolds (p_i), i.e. PNDP-manifolds that "emerge" as points (i.e., zero dimensional), and following what was said in the previous *Proposition*, we obtain that this latter can be considered as singletons that make up the spatial part of a discrete D-brane.

**Conclusions and Discussion**

We have shown that physically admissible manifolds have the same algebraic K-theories. When we know the K-theory for ordinary Quantum Field Theory, we can construct the K-theories for some more generalized Quantum Field Theories. Examples of such Quantum Field Theories is Group Field Theory. We have observed that the dynamics of strings attached on D-branes can be formulated in terms of cofibration sequences. D-branes fix a dynamical structure which the strings attached on it have.

Moreover, we have found from our K-theory condition applied to Chan-Paton bundles that only subgroups of $U(N)$ can arise as gauge groups in the low energy limit from a String theory. This can give insights on how unified field theories on energy scales below the Planck scale may look like. Further research in the role of D-branes for grand unification theories must be done to propose novel Grand Unification Theories.

In conclusion we have also shown how a special kind of Einstein sequential warped-product manifolds (called PNDP-manifolds), which use differential geometry in combination with derived geometry, can, according to some speculative interpretations, constitute the spatial part of a discrete D-brane.

**Funding**

The work was partially funded by Slovak Grant Agency for Science VEGA under the grant number VEGA 2/0076/23.